\documentclass[reqno,10pt]{amsart}
\usepackage{xspace}
\usepackage[bookmarksnumbered,colorlinks]{hyperref}
\usepackage{graphics}
\newcommand{\bt}{\begin{theorem}}                              
\newcommand{\et}{\end{theorem}}                                
\newcommand{\bd}{\begin{definition}}                           
\newcommand{\ed}{\end{definition}}                             
\newcommand{\bl}{\begin{lemma}}                                
\newcommand{\el}{\end{lemma}}                                  
\newcommand{\bpr}{\begin{proposition}}                  
\newcommand{\epr}{\end{proposition}}                    
\newcommand{\bere}{\begin{remark}}                      
\newcommand{\ere}{\end{remark}}                         
\newcommand{\beq}{\begin{equation}}
\newcommand{\eeq}{\end{equation}}
\def\dimM {r}
\def\bal#1\eal{\begin{align}#1\end{align}}              
\def\baln#1\ealn{\begin{align*}#1\end{align*}}          
\def\bml#1\eml{\begin{multline}#1\end{multline}}        
\def\bmln#1\emln{\begin{multline*}#1\end{multline*}}  
\def\bga#1\ega{\begin{gather}#1\end{gather}}
\def\bgan#1\egan{\begin{gather*}#1\end{gather*}}
\newcommand{\de}{\mathrm{d}}                        
\DeclareMathOperator{\cat}{cat}                     
\newcommand{\N}{\ensuremath{\mathbb{N}}\xspace}     
\newcommand{\R}{\ensuremath{\mathbb{R}}\xspace}     
\newcommand{\eps}{\varepsilon}                      
\newcommand{\To}{\longrightarrow}

\newcommand{\inte}{\int_0^1\!\!}
\newcommand{\tig}{\tilde{g}_0}

\newtheorem{theorem}{Theorem}[section]

\newtheorem{lemma}[theorem]{Lemma}
\newtheorem{proposition}[theorem]{Proposition}
\theoremstyle{definition}
\newtheorem{definition}[theorem]{Definition}
\theoremstyle{remark}
\newtheorem{remark}[theorem]{Remark}
\newcommand{\dist}{\ensuremath{\mathrm{dist}}\xspace}
\newcommand{\distg}{\ensuremath{\mathrm{dist}^h}\xspace}
\newcommand{\disto}{\ensuremath{\mathrm{dist}_0}\xspace}

\newcommand{\nablag}{\ensuremath{\nabla^h}\xspace}

\title[Energy functional on  Finsler manifolds]%
{On the energy functional on Finsler manifolds and applications to
stationary spacetimes}

\author[E. Caponio]{Erasmo Caponio}
\address{Dipartimento di Matematica, Politecnico di Bari, Via Orabona 4,
70125, Bari, Italy}
\email{caponio@poliba.it}
\thanks{EC supported by M.I.U.R. Research project PRIN07 "Metodi Variazionali e  
Topologici nello Studio di Fenomeni Nonlineari"
}

\author[M. A. Javaloyes]{Miguel \'Angel Javaloyes}
\address{Departamento de Geometr\'{\i}a y Topolog\'{\i}a.
  Facultad de Ciencias, Universidad de Granada.
 Campus Fuentenueva s/n, 18071 Granada, Spain}
\email{ma.javaloyes@gmail.com}
\curraddr{Departamento de Matem\'aticas, Facultad de Matem\'aticas, Universidad de Murcia,
Campus Universitario de Espinardo, 30100 Murcia, Spain}
\thanks{MAJ is partially supported by Regional J. Andaluc\'{\i}a Grant P09-FQM-4496, MICINN
project MTM2009-10418, and Fundaci\'on S\'eneca project 04540/GERM/06. This research is a result of the
activity developed within the framework of the Programme in Support of Excellence Groups of the Regi\'on de Murcia, Spain, by Fundaci\'on S\'eneca, Regional Agency for Science and Technology (Regional
Plan for Science and Technology 2007–2010).}

\author[A. Masiello]{Antonio Masiello}
\address{Dipartimento di Matematica,
Politecnico di Bari, Via Orabona 4,
70125, Bari, Italy}
\email{masiello@poliba.it}
\thanks{AM supported by M.I.U.R. Research project PRIN07 "Metodi Variazionali e  
Topologici nello Studio di Fenomeni Nonlineari"}

\subjclass[2000]{53C22, 53C50, 53C60, 58E05}

\keywords{non-reversible Finsler metrics, geodesics, stationary
Lorentzian manifolds, light rays, Fermat principle}

\begin{document}
\begin{abstract}
In this paper we first study some global properties of the energy
functional on a non-reversible Finsler manifold. In particular we
present a fully detailed proof of the Palais--Smale condition under
the completeness of the Finsler metric. Moreover we define a Finsler
metric of Randers type, which we call Fermat metric, associated to a
conformally standard stationary spacetime. We shall study the
influence of the Fermat metric  on the causal properties of the
spacetime, mainly the global hyperbolicity. Moreover we study the
relations between the energy functional of the Fermat metric and the
Fermat principle for the light rays in the spacetime.
This allows us to obtain existence and multiplicity results for light
rays, using the Finsler theory. Finally the case of timelike
geodesics with fixed energy is considered.
\end{abstract}

\maketitle

\section{Introduction}\label{sec:intro}
In the recent years there has been an increasing interest in the
study of Finsler Geometry, both from the theoretical point of view
and for the applications to many fields of Physics. We mention the
study of the multiplicity of geometrically distinct closed
geodesics, which presents different features with respect to Riemannian
Geometry as shown by the Katok's example (see \cite{Katok73}) and the study of the Zermelo
navigation problem which has led to a classification of
Randers metrics with constant flag curvature, see \cite{BaRoSh04}.
Finsler Geometry has also found many applications to applied
sciences as Biology, Classical and Quantum Optics, Relativity and
Quantum Gravity. We refer to the monographs \cite{AnInMa93},\cite{Asanov85} and to the more recent papers \cite{Duval}, \cite{GiLiSi06}, \cite{Perlic06}.

We recall some basic facts about Finsler manifolds  and we refer to
\cite{BaChSh00} for any further information. Let
$M$ be a smooth, real, paracompact manifold  of finite dimension. A
Finsler structure  on  $M$ is a function $F\colon TM\to[0,+\infty)$
which is continuous on $TM$, $C^{\infty}$ on $TM\setminus 0$,
vanishing only on the zero section, fiberwise positively
homogeneous of degree one, i.e.   $F(x,\lambda y)=\lambda F(x,y)$,
for all $x\in M$,   $y\in T_x M$ and $\lambda>0$, and which has
fiberwise strictly convex square  i.e. the matrix
\beq\label{fundamentaltensor}
g_{ij}(x,y)=\left[\frac{1}{2}\frac{\partial^2 (F^2)}{\partial
y^i\partial y^j}(x,y)\right] 
\eeq 
is positive definite for any
$(x,y)\in TM\setminus 0$. The tensor \beq\label{fundamentaltensor2}
g=g_{ij}(x,y)\de x^i\otimes\de x^j \eeq 
(here and throughout the paper we adopt the Einstein summation convention) 
is the so called {\em
fundamental tensor} of the Finsler manifold $(M,F)$; it is a
symmetric section of the tensor bundle $\pi^*(T^*M)\otimes
\pi^*(T^*M)$, where $\pi^*(T^*M)$ is the dual of the pulled-back
tangent bundle $\pi^*TM$ over $TM\setminus 0$ ($\pi$ is the
projection $TM\to M$). 
\bere\label{fquadro} We stress that, by
homogeneity, $F^2$ is $C^1$ on $TM$ and it reduces to the square  of
the norm of a Riemannian metric  if and only its  second order fiber
derivatives are continuous up to the zero section (see \cite{Warner65}). \ere 
\bere\label{nonrevers} Since $F$ is only
positive homogeneous of degree $1$,  we have that, in general,
$F(x, y)\neq F(x,- y)$. If for all $(x,y)\in TM$, $F(x,{y})=F(x,-{y})$, the
Finsler metric $F$ is said {\em reversible}. The number
$\lambda(x)=\max_{v\in T_xM}\{F(x,-{y})\ |\  F(x,{y})=1\}$ (see
\cite{Radema04}) gives a measure of  non reversibility for a Finsler
metric. \ere 
The components $g_{ij}$ of the fundamental tensor define the {\em formal Christoffel symbols} $\gamma^i_{jk}$,
\[\gamma^i_{jk}(x,y):=\frac{1}{2}g^{is}\left(\frac{\partial g_{sj}}{\partial x^k}-\frac{\partial g_{jk}}{\partial x^s}+
\frac{\partial g_{ks}}{\partial x^j}\right),\] 
and the {\em Cartan
tensor} \beq\label{cartan} A_{ijk}(x,y):=\frac{F}{2}\frac{\partial
g_{ij}}{\partial y^k}=\frac{F}{4}\frac{\partial^3(F^2)}{\partial
y_i\partial y_j\partial y_k}, \eeq for all $(x,y)\in TM\setminus 0$.
From Equation \eqref{cartan}, we see that the components $A_{ijk}(x,y)$
are totally symmetric in $(i,j,k)$ and they are positively
homogeneous of degree $0$ in the $y$ variable.

The {\em Chern connection} $\nabla$ is the unique linear connection on $\pi^* TM$ whose connection $1$-forms
$\omega^i_j$ satisfy the following structural equations:
\begin{align}
&\de x^j\wedge \omega^i_j=0&&\text{\em torsion free,}\label{torsionfree}\\
&\de g_{ij}-g_{kj}\omega^k_i-g_{ik}\omega^k_i=\frac{2}{F}A_{ijs}\delta y^s&&\text{\em almost $g$-compatibility,}\label{almostgcomp}
\end{align}
where $\delta y^s$ are the $1$-forms on $\pi^* TM$ given as
$\delta y^s:=\de y^s + N^s_j\de x^j$, and
\[N^i_j(x,y):=\gamma^i_{jk}y^k-\frac{1}{F}A^i_{jk}\gamma^k_{rs}y^ry^s\]
are the coefficients of the so called {\em nonlinear connection} on $TM\setminus 0$.
The components of the Chern connection are given by:
\beq\label{Cherncomponents}
\Gamma^i_{jk}(x,y)=\gamma^i_{jk}-\frac{g^{il}}{F}\left(A_{ljs}N^s_k-A_{jks}N^s_i+A_{kls}N^s_j\right).
\eeq
Clearly  $\Gamma^i_{jk}(x,y)$ are defined on $TM\setminus 0$  and they are positively homogeneous of degree $0$ with respect to $y$.

Let $\gamma=\gamma(s)$ be a smooth regular curve on $M$, with velocity field $T=\dot\gamma$, and $W$ be a   smooth vector field along $\gamma$.
The Chern connection defines two different covariant derivatives $D_TW$ along $\gamma$:
\bal
&D_{T}W=\left.\left(\frac{\de W^i}{\de t}+W^jT^k\Gamma^i_{jk}(\gamma,T)\right)\frac{\partial}{\partial x^i}\right\vert_{\gamma(t)}&&\text{\em with reference vector $T$,}\label{referenceT}\\
&D_{T}W=\left.\left(\frac{\de W^i}{\de t}+W^jT^k\Gamma^i_{jk}(\gamma,W)\right)\frac{\partial}{\partial x^i}\right\vert_{\gamma(t)}&&\text{\em with reference vector $W$.}\nonumber
\eal
A geodesic of the Finsler manifold $(M,F)$ is a smooth regular curve $\gamma$ satisfying
the equation
\[
D_{T}\left(\frac{T}{F(\gamma,T)}\right)=0,
\]
with reference vector $T=\dot\gamma$. A curve $\gamma=\gamma(s)$ is
said to have {\em constant speed} if $F(\gamma(s),\dot \gamma (s))$
is constant along $\gamma$. Constant speed geodesics satisfy the
equation \beq\label{geoconstspeed} D_T T=0, \eeq with reference
vector $T=\dot\gamma$. The length of a piecewise smooth  curve
$\gamma\colon [a,b]\subset\R \to M$ with respect to the Finsler
structure $F$  is defined by
\[L(\gamma)=\int_a^b\!\!
F(\gamma(s),\dot\gamma(s))\de s.\] Thus the distance  between two
arbitrary points $p,\ q\in M$ is given by \beq\label{finsleriandist}
\dist(p,q)= \inf_{\gamma\in C(p,q)}L(\gamma), \eeq where $C(p,q)$ is
the set of all piecewise smooth curves $\gamma\colon[a,b]\to M$ with
$\gamma(a)=p$ and $\gamma(b)=q$. The distance function
\eqref{finsleriandist} is nonnegative and satisfies the triangle
inequality, but it is not symmetric as $F$ is non-reversible. Thus
one has  to distinguish the order  of   a pair  of points in $M$
when speaking about  distance. As a consequence, one is naturally
led to the notions of forward and backward metric balls, spheres,
Cauchy sequences and completeness (see \cite[\S 6.2]{BaChSh00}). For
instance: the {\em forward metric ball} $ B^+_r(p)$ (resp. {\em
backward} $ B^-_r(p)$) of center $p\in M$ and radius $r\geq 0$ is
given by all the points $x\in M$ such that $\dist(p,x)<r$ (resp.
$\dist (x, p)< r$); a sequence $\{x_n\}\subset M$ is called {\em
forward} (resp. {\em backward}) {\em Cauchy sequence} if for all
$\eps>0$ there exists $\nu\in\N$ such that, for all $\nu\leq i\leq
j$, $\dist(x_i,x_j)\leq \eps$ (resp. $\dist(x_j,x_i)\leq \eps$);
$(M,F)$ is {\em forward} (resp. {\em backward}) {\em complete} if
all forward (resp. backward) Cauchy sequences are convergent;
$(M,F)$   is said {\em forward} (resp. {\em backward}) {\em geodesically complete} if every geodesic $\gamma\colon[a,b)\to M$ (resp.
$\gamma\colon (b, a]\to M$) can be extended to a geodesic defined on
the interval $[a,+\infty)$ (resp. $(-\infty, a]$). What is relevant
here is that the topologies generated by the forward and the
backward metric balls coincide with the underlying manifold
topology; moreover a suitable version of Hopf-Rinow theorem holds
(see \cite[Theorem 6.6.1]{BaChSh00}) stating the equivalence of the
notions  of forward (or backward) completeness and the compactness of
closed and forward (or backward) bounded subsets of $M$ and implying the
existence of a geodesic connecting any pair of points in $M$
and minimizing the Finslerian distance. It is worth recalling that
the two notions of completeness are not equivalent (see  for example
\cite[\S 12.6.D]{BaChSh00}).

As in Riemannian Geometry, geodesics on a Finsler manifold $(M,F)$
satisfy a variational principle. First of all a curve is a geodesic
for the Finsler metric $F$ if and only if it minimizes the length
between two sufficiently close points on the curve, see
\cite{BaChSh00}. Moreover a smooth curve $x\colon
[a,b]\rightarrow M$ is a constant speed geodesic if and only if it
is a stationary point of the {\em energy functional}
\[
J(x) = \frac 12\int_a^b F^2(x,\dot x){\rm d}s
\]
on the space of sufficiently smooth curves on the manifold $M$
joining the points $x(a)$ and $x(b)$ (for more general
boundary conditions see Section 2).

In this paper we shall study the main properties of the energy
functional in the infinite dimensional setting of the
Sobolev-Hilbert manifold of $H^1$ curves satisfying very general
boundary conditions, containing the classical two points and
periodic boundary conditions. In particular we shall present a fully
detailed proof of the Palais--Smale condition for the energy
functional. In the second part of the paper we present a new
application of Finsler Geometry to General Relativity. In the class
of conformally standard stationary spacetimes we define a Finsler
metric of Randers type, which we call {\em Fermat metric}. The
choice of this definition is due to the fact that this metric is
strictly related to the Fermat principle of light rays in this class
of spacetimes. We shall also show that the causal structure of a
conformally stationary spacetime is influenced by the global
properties of the Fermat metric. In particular the global
hyperbolicity of the metric is strictly related to the completeness
of the Fermat metric. Moreover the equivalence between the Fermat
principle of light rays and the geodesic problem for the Fermat
metrics allows one to obtain multiplicity results for light rays as an
application of the analogous results in the Finsler setting. These
results allow  us to obtain a mathematical model of the gravitational
lens effect on conformally stationary spacetimes. Finally 
analogous results for timelike geodesics on a stationary spacetime are
presented.

\section{The energy functional in Finsler Geometry}\label{mercuri}
In this section we shall study the energy functional  of a Finsler
manifold $(M,F)$ in the infinite dimensional setting of
Hilbert-Sobolev manifolds. We recall that the infinite dimensional
setting for the energy functional and the variational theory for
geodesics on a Riemannian manifold was introduced by R. Palais in
the paper \cite{Palais63} and extended by F.Mercuri to Finsler
manifolds in the paper \cite{Mercur77}. Here we shall prove in all
the details that the critical points of the energy functional,
defined on a manifold of curves satisfying boundary
conditions \eqref{orto}, are smooth and they are exactly the geodesics satisfying (\ref{orto}) and parametrized with
constant speed. 

Let $(M,F)$ be a  forward or  backward complete Finsler manifold and
let us endow $M$ with any complete Riemannian metric $h$. Let $N$ be
a smooth submanifold of $M\times M$. We consider the collection
$\Lambda_N(M)$ of curves $x$ on $M$ parameterized on the interval
$[0,1]$ with endpoints  $\big(x(0),x(1)\big)$  belonging to $N$ and
having $H^1$ regularity, that is, $x$ is absolutely continuous and
the integral $\inte h(x)[\dot x,\dot x]\de s$ is finite.  It is well
known that $\Lambda_N(M)$ is a Hilbert manifold modeled on any of
the equivalent  Hilbert  spaces of all the $H^1$ sections, with endpoints in
$TN$, of the pulled back bundle $x^*TM$, $x$  any  regular curve in
$\Lambda_N(M)$, \cite[Proposition 2.4.1]{Klinge82}. In fact, the
scalar product is given by
\begin{equation}\label{hilbertproduct}
\langle X,Y\rangle_1=\inte h(x)[X,Y]\de s+\inte h(x)[\nabla^h_x X,\nabla^h_xY]\de s,
\end{equation}
for every  $H^1$ sections, $X$ and $Y$ of $x^*TM$,
$\nablag_x$ being the covariant derivative along $x$ associated to the Levi-Civita connection of the metric $h$.

Let us denote the function $F^2$ by $G$ and let us consider the {\em
energy functional} \beq\label{J} J(x)=\frac{1}{2}\inte G(x,\dot
x)\de s \eeq of the  Finsler manifold $(M,F)$,
defined on the manifold $\Lambda_N(M)$. The functional $J$ is
$C^{2-}$ on $\Lambda_N(M)$, i.e. it is $C^1$ with locally Lipschitz
differential (see \cite[Theorem 4.1]{Mercur77}).\footnote{Though in
that paper $J$ is defined on $\Lambda_{\triangle}(M)$, where
$\triangle$ is the diagonal in $M\times M$, the $C^{2-}$ regularity of $J$ on $\Lambda_N(M)$ can be carried out along the same lines.} A  {\em critical point} $\gamma$
of $J$ is a curve $\gamma\in\Lambda_N(M)$ such that $\de
J(\gamma)=0$.

We first study the regularity of critical points for $J$. We shall
show in all the details that, in spite of the loss of regularity of
the Lagrangian function $G$ on the zero section, the $H^1$--critical points of $J$ are
smooth curves.

We start by  computing the differential of $J$ on $\Lambda_N (M)$  to show that a  non constant critical point is a geodesic
satisfying the boundary conditions
\beq\label{orto}
g\big(\gamma(0),\dot\gamma(0)\big)[V,\dot \gamma(0)]=g\big(\gamma(1),\dot\gamma(1)\big)[W,\dot \gamma(1)],
\eeq
where $g$ is the fundamental tensor of the Finsler metric $F$ defined in \eqref{fundamentaltensor} and  $(V,W)\in T_{(\gamma(0),\gamma(1))}N$.
\bere
Before the next Lemma, let us see what the boundary conditions \eqref{orto} become in some particular cases:
\begin{itemize}
\item[(i)] Let $\triangle$ be the diagonal in $M\times M$ and $N=\triangle$. From \eqref{fundamentaltensor} and the Euler theorem for homogeneous functions, we know that
$\partial_yG(x,y)= 2 g(x,y)[\cdot,y]$, for any $(x,y)\in TM$. Hence, from  $\gamma(0)=\gamma(1)$ and  \eqref{orto} we get
\[\partial_yG(\gamma(0),\dot\gamma(0))=\partial_yG(\gamma(0),\dot\gamma(1)).\]
Since the map $y\mapsto\partial_yG(x,y)$ is an injective map (see the proof of Proposition~\ref{regularity} below), it has to be $\dot\gamma(0)=\dot\gamma(1)$.
\item[(ii)] Let $M_0$ and $M_1$ be two submanifolds of $M$ and $N=M_0\times M_1$. In \eqref{orto} take $W=0$. Then, for any $V\in T_{\gamma(0)}M_0$ we get $g\big(\gamma(0),\dot\gamma(0)\big)[V,\dot \gamma(0)]=0$. Analogously taking $V=0$, it has to be
$g\big(\gamma(1),\dot\gamma(1)\big)[W,\dot \gamma(0)]=0$, for any
$W\in T_{\gamma(1)}M_1$.
\end{itemize}
\ere

\bl\label{derivate}
Let $\gamma\colon[0,1]\to M$ be a smooth regular curve and $\sigma\colon[0,1]\times[-\eps,\eps]\to M$, $\sigma=\sigma(t,u)$ be a smooth regular variation of $\gamma$ (i.e. $\sigma(t,0)=\gamma(t)$ for all $t\in[0,1]$) with variation vector field $U=\partial_u\sigma$. Then
\beq
\partial_u \left(g(\sigma, T)[T,T]\right)=2g(\sigma,T)[T,D_UT],\label{productrule}
\eeq
where $T=\partial_t\sigma$ and the covariant derivative $D_UT$ is with reference vector $T$ (see formula \eqref{referenceT}).
\el
\begin{proof}
From the symmetry of the $g_{ij}$ we get
\bal
\lefteqn{\partial_u \left(g_{ij}(\sigma, T)T^iT^j\right)=}&\nonumber\\
&=\partial_{x^k}\left(g_{ij}(\sigma,T)\right)U^kT^iT^j+\partial_{y^k}\left(g_{ij}(\sigma,T)\right)(\partial_u T)^kT^iT^j\nonumber\\&\quad+2g_{ij}(\sigma,T)T^i\left(\partial_u T\right)^j,\label{partialu}\\
\intertext{and using the definition of the Cartan tensor \eqref{cartan}, the right-hand side of \eqref{partialu} becomes}
&\partial_{x^k}\left(g_{ij}(\sigma,T)\right)U^kT^iT^j+\frac{2A_{ijk}(\sigma,T)}{F(\sigma,T)}(\partial_u T)^kT^iT^j+2g_{ij}(\sigma,T)T^i\left(\partial_u T\right)^j.\label{productrule1}
\eal
On the other hand, by equating the coefficients of the $1$-forms $\de x^k$ in \eqref{almostgcomp} we see that
\beq\label{fromalmostgcomp}
\partial_{x_k}g_{ij}=g_{sj}\Gamma^s_{ik}+g_{is}\Gamma^s_{jk}+\frac{2}{F}A_{ijs}N^s_k.
\eeq
Now  we recall that $G$ is positively homogeneous of degree $2$ in $y$ and consequently, from Eq. \eqref{cartan} and Euler's theorem, we get
\beq\label{eulerA}
y^iA_{ijk}(x,y)=y^jA_{ijk}(x,y)=y^kA_{ijk}(x,y)=0.
\eeq
Hence the terms
\baln
\frac{2A_{ijk}(\sigma,T)}{F(\sigma,T)}(\partial_u T)^kT^iT^j&& \text{and}&&\frac{2A_{ijk}(\sigma,T)}{F(\sigma,T)}T^iT^jN^s_k(\sigma,T)U^k,
\ealn
appearing after substituting \eqref{fromalmostgcomp} in \eqref{productrule1}, vanish. Finally,
using again the symmetry of $g_{ij}$ we obtain
\[
\partial_u \left(g_{ij}(\sigma, T)T^iT^j\right)=2g_{ij}(\sigma,T)T^i\left((\partial_uT)^j+\Gamma^j_{hk}(\sigma,T)T^hU^k\right),
\]
which is Eq. \eqref{productrule} in local coordinates.
\qed\end{proof}
Now we can prove the following.
\bpr\label{regularity}
A   curve  $\gamma\in\Lambda_N(M)$ is  a constant (non zero) speed geodesic for the Finsler manifold $(M,F)$ satisfying
\eqref{orto}
if and only if it is a (non constant) critical point of $J$.
\epr
\begin{proof}
Let $\gamma\colon[0,1]\to M$ be a   smooth curve in $\Lambda_N(M)$ and $Z\in T_{\gamma}\Lambda_N(M)$ be a smooth vector field along $\gamma$.
Choose a smooth variation $\sigma\colon[0,1]\times[-\eps,\eps]\to M$, $\sigma=\sigma(t,u)$ of $\gamma$ with variation vector field $U=\partial_u\sigma$ having endpoints in $TN$ and such that $U(t,0)=Z(t)$ for all $t\in[0,1]$. Moreover we set $T=\partial_t\sigma$.
We will cover the support of $\gamma$ by a finite number of local charts $\{(V_k,\varphi_k)\}$ of the manifold $M$ so that the variation of $J$ will be computed using the  systems of coordinates induced on $TM$.
 With abuse of notations we shall not change the symbols denoting  the points, the vectors  and the forms in such coordinate systems and we shall omit the sum symbol in the integrands.
Since $G$ is $C^1$ on $TM$ (see Remark \ref{fquadro})  and using the equality  $\partial_uT=\partial_t U$, we get
\begin{align*}
\frac{\de}{\de u}J(\sigma)&=\frac{1}{2}\!\inte\!\partial_u\big(G(\sigma, T)\big)\de t=\frac{1}{2}\!\inte\!\big(\partial_xG(\sigma,T)[U]+\partial_yG(\sigma,T)[\partial_uT]\big)\de t\\&=
\frac{1}{2}\inte \big(\partial_xG(\sigma,T)[U]+\partial_yG(\sigma,T)[\partial_tU]\big)\de t,
\end{align*}
which evaluated at $u=0$ gives \beq\label{1variation} \de
J(\gamma)[Z]=\frac{1}{2}\inte\big(\partial_x
G(\gamma,\dot\gamma)[Z]+\partial_yG(\gamma,\dot\gamma)[\dot
Z]\big)\de t. \eeq 
This equation can be extended by
density to any curve $\gamma\in\Lambda_N(M)$ and to any vector field
$Z\in T_{\gamma}\Lambda_N(M)$. Now assume that
$\gamma\in\Lambda_N(M)$ is a critical point of $J$. We are going to
show that $\gamma$ is a smooth curve. Evaluating \eqref{1variation}
on any smooth vector field $Z$ with compact support in the interval
$I_k=\gamma^{-1}(V_k)=(t_k,t_{k+1})\subset [0,1]$ we get
\beq\label{dubois}\int_{I_k}\big( H
+\partial_yG(\gamma,\dot\gamma)\big)[\dot  Z]\de t=0, \eeq where
$H=H(t)$ is the $T^*M$ valued  function
\[
H(t)=-\int_{t_k}^t\big(\partial_x G(\gamma,\dot\gamma)\big)\de s.
\]
Last integration has only a local sense, since it consists of the
integrals of the components of the covector $\partial_x
G(\gamma,\dot\gamma))$ along the curve $\gamma$. Moreover, equation
\eqref{dubois} implies that there exists a constant  covector
$W\in\mathbb{R}^n$, with $n=\dim M$, such that
\begin{equation}\label{integraloperator}H(t)+\partial_yG(\gamma(t),\dot\gamma(t))=W,
\end{equation} a. e. on $I_k$; since $H$ is continuous, the function $t\in I_k\mapsto \partial_yG(\gamma,\dot\gamma)$ is also continuous.
Now fix $x\in M$ and consider  the map $\mathcal L_x:=y\in T_x
M\setminus \{0\} \mapsto \partial_y G(x,y)\in T^*_{x}M$.  Since $G$
vanishes only on the zero section and  is positively  homogeneous of
degree $2$ in $y$, by Euler's theorem also $\partial_y G(x,y)$
is the map of constant value $0$ if and only if $y=0$. Hence $\mathcal L_x$ assumes
values in $T_x^*M\setminus\{0\}$. Being $G$ fiberwise strictly
convex on $TM\setminus 0$, $\mathcal L_x$ is locally invertible on
$T_xM\setminus\{0\}$ . Moreover, as $\mathcal L_x$ is positively
homogeneous of degree $1$, it is a proper map and therefore it is a
homeomorphism from $T_x M\setminus\{0\}$ onto $T_x^*
M\setminus\{0\}$ (see \cite[Theorem 1.7, p. 47]{AmbPro93}). Since
the inverse of a homogeneous function of degree $1$ is homogeneous
of degree $1$   and $\mathcal L_x(y)=0$ if and only if $y=0$, we
obtain that $\mathcal L_x$ is a homeomorphism from $T_x M$ onto
$T^*_x M$. Now consider the maps $\Phi\colon=(x,y)\in TM\mapsto
(x,\mathcal L_x(y))\in T^*M$ and $\Psi\colon=(x,w)\in T^*M\mapsto
(x,\mathcal L^{-1}_x(w))\in TM$. As $\partial_{yy} G$ is positive
definite on $TM\setminus 0$, from the inverse function  theorem
$\Phi$ is locally a homeomorphism on $TM\setminus 0$ and
$\Phi^{-1}=\Psi$ on $T^*M\setminus 0$. Therefore  the map $(x,w)\in
T^*M\setminus 0\mapsto \mathcal L^{-1}_x(w)$ is continuous and the
continuity extends up to the zero section. In fact  if $(x_n,w_n)\to
(\bar x, 0)$ then
\[
\mathcal L^{-1}_{x_n}(w_n)=  |w_n|\mathcal L^{-1}_{x_n}\Big(\frac{w_n}{|w_n|}\Big)\to 0,
\]
where we have identified a neighborhood of $(\bar x,0)$ on $TM$ or $T^*M$ with an open set of $\R^n\times\R^n$, and we have used the continuity of the map $\mathcal L^{-1}$ on $T^*M\setminus 0$.
 Thus, we can state that the function $t\in I_k\mapsto \mathcal L^{-1}_{\gamma(t)}\circ\mathcal L_{\gamma(t)}(\dot\gamma(t))=\mathcal L^{-1}_{\gamma(t)}\big(\partial_yG(\gamma(t),\dot\gamma(t))\big)=\dot\gamma(t)$ is continuous
and  $\gamma$ is a $C^1$ curve. From  \eqref{integraloperator}, we get that
$\gamma$  satisfies the following equation a. e. on $I_k$
\begin{equation}\label{lagrange}
\frac{\de}{\de t}\partial_yG(\gamma,\dot\gamma)=\partial_xG(\gamma,\dot\gamma).
\end{equation}
Hence we deduce that $\frac{\de}{\de t}\partial_yG(\gamma,\dot\gamma)$ is continuous on $I_k$. This  information and the fact that $G$ is fiberwise strictly convex imply  that $\gamma$ is actually twice differentiable  on every point $t$ where $\dot\gamma(t)\neq 0$ (see for instance \cite[Proposition 4.2]{BuGiHi98}).
Now assume that $\gamma$ is not a constant curve and let $A_k\subset I_k$ be the open subset of the points $t\in I_k$ where $\dot\gamma(t)\neq 0$. From \eqref{lagrange} we see that  the energy $E(\gamma):=\partial_y G(\gamma,\dot\gamma)[\dot\gamma]-G(\gamma,\dot\gamma)$ is constant on every connected component of $A_k$. Since $G$ is positively homogeneous of degree $2$, from the Euler's theorem,
we have  that $E(\gamma)=G(\gamma,\dot\gamma)$. Recalling that
$F$ is zero only on the zero section and that the function $t\in I_k\mapsto G(\gamma(t),\dot\gamma(t))$  is continuous, we conclude
that $G(\gamma(t),\dot\gamma(t))$ is constant (non zero) on every $I_k$.  As we can enlarge all the intervals $I_k$, except the last one, a small $\epsilon$, all the constants have to be the same  and therefore $\gamma$ is a smooth regular curve.

Now let  $Z\in T_{\gamma}\Lambda_N(M)$ be a smooth vector field along $\gamma$ and  let  $\sigma\colon[0,1]\times[-\eps,\eps]\to M$, $\sigma=\sigma(t,u)$ be a smooth regular variation of $\gamma$ with variation vector field $U=\partial_u\sigma$ having endpoints in $TN$ and such that $U(t,0)=Z(t)$ for all $t\in[0,1]$.
Since $G(x,y)=g(x,y)[y,y]$ for any $(x,y)\in TM\setminus 0$,
from \eqref{productrule} we get
\beq\label{new1variation}
\frac{\de}{\de u}J(\sigma)=\frac{1}{2}\inte\partial_u\left(g(\sigma,T)[T,T]\right)\de t=\inte g(\sigma,T)[T, D_UT]\de t,
\eeq
where $T=\partial_t\sigma$. On the other hand, as the variation $\sigma$ is smooth, it holds $D_UT=D_TU$ both with reference vector $T$ and hence, using this equality in \eqref{new1variation} and evaluating at $u=0$, we obtain
\beq\label{differential}
\de J(\gamma)[Z]=\inte g(\gamma,\dot\gamma)[\dot\gamma,D_{\dot\gamma}Z]\de t,
\eeq
where $D_{\dot\gamma}Z$ has reference vector $\dot\gamma$. Moreover,
arguing as in the proof of Lemma~\ref{derivate}, one gets
\[\frac{\de}{\de t}\left(g(\gamma,\dot\gamma)[\dot\gamma,Z]\right)=g(\gamma,\dot\gamma)[D_{\dot\gamma}\dot\gamma,Z]+g(\gamma,\dot\gamma)[\dot\gamma,D_{\dot\gamma}Z],\]
which, when applied to \eqref{differential}, gives us
\bal
0=\de J(\gamma)[Z]&=-\inte g(\gamma,\dot\gamma)[D_{\dot\gamma}\dot\gamma,Z]\de t\nonumber\\
&\quad+g(\gamma(1),\dot\gamma(1))[\dot\gamma(1),Z(1)]-
g(\gamma(0),\dot\gamma(0))[\dot\gamma(0),Z(0)].\label{postdifferential}
\eal
Finally, by choosing  an endpoints vanishing vector field $Z$ we
see that $\gamma$ has to satisfy the equation $D_{\dot\gamma}\dot\gamma=0.$
Thus, $\gamma$ is a constant speed geodesic satisfying
the boundary conditions  \eqref{orto}.

For the converse, we observe that  if $\gamma$ is a constant
non-zero speed geodesic satisfying the boundary conditions
\eqref{orto}, then \eqref{postdifferential} holds and  hence  $\gamma$ is a critical point of $J$.
\qed\end{proof}

\section{On the Palais-Smale condition for the energy functional}\label{PS}
We prove now that the energy functional $J$ satisfies the
Palais--Smale condition. We recall that a functional $J$ defined on
a Banach manifold $X$ satisfies the Palais-Smale condition if every
sequence $\{x_n\}_{n\in N}$ such that $\{J(x_n)\}_{n\in \N}$ is
bounded and $\|\de J(x_n)\|\to 0$ contains  a convergent
subsequence.

The lack of regularity of the function $G = F^2$ on the zero section
gives rise to some problems, for instance in the application of the
mean value theorem, which do not occur in the proof of the
Palais-Smale condition for the energy functional of a Riemannian
manifold (see for instance \cite{Klinge82}). In the paper \cite{Mercur77} such problems are
circumvented by using a sketched approximation argument. Here we
give a fully detailed proof of the Palais--Smale condition. By a
{\em localization argument} we will work on an open subset of
$\R^n$. This allows us to reduce the technical aspects of the proof.

 \bt\label{ps} Let $(M,F)$ be
forward (resp. backward) complete and $N$ be  a closed submanifold on
$M\times M$ such that the first projection (resp. the second
projection) of $N$ to $M$ is compact, then $J$ satisfies the
Palais-Smale condition on $\Lambda_N(M)$. \et
\begin{proof}
We   prove the theorem in the forward complete case, being the backward one
analogous.
Let $\{x_n\}_{n\in \N}$ be a sequence contained in $\Lambda_N(M)$ on which $J$
is bounded.
Under the assumptions of Theorem~\ref{ps}, the manifold
$\Lambda_N(M)$ is a closed submanifold of the complete Hilbert manifold
$\Lambda(M)$ (see \cite[Theorem 2.4.7]{Klinge82}), which is   the collection of all the $H^1$ curves in $M$
parameterized on the interval $[0,1]$  with scalar product as in Eq. \eqref{hilbertproduct}.
The  differentiable manifold structure on $\Lambda(M)$ is given by the charts
$\{(U_{\omega},\exp^{-1}_{\omega})\}_{\omega\in C^{\infty}(M)}$, where $\exp^{-1}_{\omega}$ is the
inverse of the map $\exp_{\omega}(\xi)=\exp_{\omega(t)}\xi(t)$, for all $\xi\in H^1(O_{\omega})$, being $\exp\colon TM\to M$ the exponential map of the Riemannian manifold $(M,h)$ and $O_{\omega}$ a neighborhood of the zero
section in ${\omega}^*TM$ (see \cite[Theorem 2.3.12]{Klinge82}).

First we prove that $\{x_n\}$ converges uniformly.
Pick a point  $\bar p\in \mathrm{p}_1(N)$, where $\mathrm{p}_1$
is the first projection of $M\times M$. We  evaluate the distance
\baln
\dist\big(\bar p, x_n(s)\big)&\leq \dist\big(\bar p,
x_n(0)\big)+\dist\big(x_n(0), x_n(s)\big)\\
                                                &\leq  \dist\big(\bar p,
x_n(0)\big)+\inte F(x_n,\dot x_n)\de s,
\ealn
for all $s\in[0,1]$, $n\in\N$. Since $\mathrm{p}_1(N)$ is compact, there
exists a constant $K$ such that
$\dist\big(\bar p, x_n(0)\big)\leq K$. By the H\"older inequality we get
\[\dist\big(\bar p, x_n(s)\big)\leq K+\left(\inte G(x_n,\dot x_n)\de
s\right)^{\frac{1}{2}}\leq K_1.\]
Then by the Finslerian Hopf-Rinow theorem the supports of the curves $x_n$ are
contained in a compact subset $C$ of $M$.  Hence there exist two  positive
constants $c_1$ and $c_2$ such that
$c_1|y|^2\leq G(x,{y})\leq c_2|{y}|^2$, for every  $x\in C$ and for every ${y}\in
T_x M$.  Here $|\cdot|$ is the norm associated to
the metric $h$. Moreover,  let $\dist_h$ be the distance associated to the Riemannian metric $h$, then  using the last inequality and again the H\"older's one, we get

\begin{align*}\dist_h(x_n(s_1),x_n(s_2))\leq &\int_{s_1}^{s_2}|\dot x_n|\de s\leq \sqrt{s_2-s_1}\left(\int_0^1|\dot x_n|^2\de s\right)^{\frac 12}\\
\leq & \frac{1}{c_1} \sqrt{s_2-s_1}\left(\int_0^1 G(x_n,\dot x_n)\de s\right)^{\frac 12}\leq K_2\sqrt{s_2-s_1},
\end{align*}
 with $s_1<s_2$ in $[0,1]$ and $K_2>0$.  Hence $\{x_n(t)\}$ is relatively compact for every $t\in[0,1]$  and uniformly   H\"older. Therefore
we can use the symmetric distance induced by $h$ and the Ascoli-Arzel\`a
theorem to obtain a subsequence,   which will be denoted again by $\{x_n\}$,
converging uniformly
to a $C^0$ curve $\bar x$ parameterized in $[0,1]$ and having endpoints in $N$.

Now we introduce the {\it localization argument} as in Appendix A.1 of \cite{AbbFig07}. Given any 
$\eta>0$ small enough we have that the subset $\mathcal{C}=\{{\rm exp}_{\bar{x}(s)}v:s\in[0,1];v\in 
\bar{B}({0},\eta)\subset T_{\bar{x}(s)}M\}$ is compact. Let $\mu(p)$ be the injectivity radius of $p$ in $(M,h)$ and $\rho=\inf\{\mu(p):p\in\mathcal{C}\}$. As the injectivity radius is continuous (see \cite[Proposition 8.4.1]{BaChSh00}), $\rho>0$ and we can choose a $C^\infty$ curve $\omega$  in such a way that $\|\bar{x}-\omega\|_{\infty}<\min\{\eta,\rho/2\}$.  Let $[0,1]\ni t\rightarrow{\mathbf E}(t)=(E_1(t),\ldots,E_r(t))$ be a parallel orthonormal frame along $\omega$, with $r=\dim M$, $P_t:\R^r\rightarrow T_{\omega(t)}M$  defined as $P_t(v_1,\ldots,v_r)=v_1 E(t)+\ldots +v_r E_r(t)$ and consider the  Euclidean open ball of radius $\rho$,  which we name $U$, and the map
$\varphi(t,v)={\rm exp}_{\omega(t)}P_t(v)$. As $\rho$ is smaller than the injectivity radius of $\omega(t)$,  the map $\varphi_t:U\rightarrow M$, defined as $\varphi_t(v)=\varphi(t,v)$,  is locally invertible and injective with invertible differential $\de \varphi_t(v)$, for every $t\in[0,1]$ and  $v\in U$.
By taking a smaller open in $U$ that contains the closed ball of radius $\rho/2$ and it is contained in the closed ball of radius $2\rho/3$,
we can assume that all the continuous functions involved in the rest of the proof  are uniformly bounded in 
$[0,1]\times U$ or in $\cup_{t\in[0,1]}\{t\}\times\varphi (\{t\}\times U)$, as for example the norms of $\de\varphi(t,v)$ and $\de \phi(t,x)$, where $\phi(t,x)=\varphi_t^{-1}(x)$.
Let ${\mathcal O}_{\omega}$ be a neighborhood of $\omega$ in $H^1([0,1],M)$  such that the map
\[\varphi_*^{-1}: {\mathcal O}_{\omega}\rightarrow H^1([0,1],U),\]
defined as $\varphi_*^{-1}(x)(t)=\varphi_t^{-1}(x(t))$ is the map of a coordinate system centered at $\omega$. Observe that the inverse of $\varphi_*^{-1}$ is the map $\varphi_*$, defined by $\varphi_*(\xi)(t)=\varphi(t, \xi(t))$. Clearly if $n$ is big enough, $x_n\in \varphi_*(H^1([0,1],U))$ and we call $\xi_n=\varphi^{-1}_*(x_n)$. Hence, proving the strong convergence of $\{x_n\}$ is equivalent to proving the strong convergence of $\{\xi_n\}$ in $H^1([0,1],U)$.

Now
consider the orthogonal splitting  \[H^1\big([0,1],\R^{\dimM}\big)=H^1_0\big([0,1],\R^{\dimM}\big)\oplus V,\]
where  $V$
is the  vector space of dimension $2r$, defined as $V=\{\zeta\in
C^{\infty}([0,1],\R^{\dimM})\ |\zeta''-\zeta=0\}$.
So, if $n\in\N$ is big enough there exist $\xi_n^0\in H^1_0([0,1],U)$ and
$\zeta_n\in V$ such that $\xi_n=\xi_n^0+\zeta_n$.
Considering $J$ as defined on $H^1([0,1],M)$, we have:
\bal
\lefteqn{\de(J\circ\varphi_*)(\xi_n)[\xi_n-\xi_m]}&\nonumber\\
&=\de(J\circ\varphi_*)(\xi_n)[\xi^0_n-\xi^0_m]+\de(J\circ\varphi_*)(\xi_n)
[\zeta_n-\zeta_m]\nonumber\\
&=\de J(x_n)\big[\de \varphi_*(\xi_n)[\xi^0_n-\xi^0_m]\big]+
\de(J\circ\varphi_*)(\xi_n)
[\zeta_n-\zeta_m]\To 0,\label{byps}
\eal
as $n,m\to \infty$. Indeed, the first term on the right-hand side of \eqref{byps} goes to $0$ as $n,m\to \infty$ since $\{x_n\}$ is a Palais-Smale
sequence  for $J$ on $\Lambda_N(M)$, the norms of the operators $\de \varphi_*(\xi_n)$ are uniformly bounded and $\{\xi_n\}$ is a bounded sequence in $H^1([0,1],U)$ (and hence also $\{\xi_n^0\}$ is a bounded sequence in
$H^1_0([0,1],U)$).  The fact that $\{\xi_n\}$ is a bounded sequence in
$H^1([0,1],U)$ follows from the inequality
\begin{multline}\label{xil2bound}
\inte|\dot \xi_n|^2\de s=\inte|\de \phi(s,x_n)[(1,\dot x_n)]|^2\de s\\
\leq K_3\inte (1+h(x_n)[\dot x_n,\dot x_n])\de s\leq
K_3 +K_4J(x_n)<K_5<+\infty,
\end{multline}
where $\phi(s,x)=\varphi^{-1}_s(x)$, for each $s\in[0,1]$ and $x\in \varphi_s(U)$, and $K_3,K_4,K_5$ are positive constants.
The second term on the right-hand side of \eqref{byps} goes also to $0$ as it can be easily seen observing that $\{\zeta_n\}$ is a converging sequence in the $C^1$ norm (this follows from the
$C^0$ convergence of $\{\xi_n\}$ and the smooth dependence of the solutions of the differential equation defining $V$ on boundary data) and using \eqref{carta} below,
with $\zeta_n-\zeta_m$ in place of $\xi_n-\xi_m$, together with the fact that $\{\xi_n\}$ is bounded in $H^1([0,1],U)$.

To complete the proof, we have to show that the sequence of curves $\{\xi_n\}$ is Cauchy in the $H^1$ norm. Notice that $\tilde{J}=J\circ\varphi_*$ is given by
$\tilde{J}(\xi)=\frac 12\inte \tilde{G}_s(\xi,\dot\xi)\de s$
for $\xi\in H^1([0,1],U)$, where $\tilde{G}_s:U\times\R^r\rightarrow \R$ is defined as
\[\tilde{G}_s(x,y)=
 G\big(\varphi(s,x),\de \varphi(s,x)[(1,y)]\big).\]
By \eqref{byps} we have
\bml
\de \tilde{J}(\xi_n)[\xi_n-\xi_m]\\
=\frac{1}{2}\inte \partial_x\tilde{G}_s(\xi_n,\dot  \xi_n)[\xi_n-\xi_m]\de s
+\frac{1}{2}\inte\partial_y \tilde{G}_s(\xi_n,\dot \xi_n)[\dot \xi_n-\dot \xi_m]\de s\to 0,\label{carta}
\eml
as $m$ and $n$ go to $\infty$. Now consider the first integral in \eqref{carta}.
We observe that,  with the same abuse of notation as in the proof of Proposition~\ref{regularity},
\begin{multline*}
\partial_x \tilde{G}_s(x,y)[\cdot]\!=\!\partial_x G(\varphi(s,x),\de \varphi(s,x)[(1,y)])[\de\varphi_s(x)[\cdot]]\\
+\partial_y G(\varphi(s,x),\de\varphi(s,x)[(1,y)])\big[\partial^2_{sx}\varphi(s,x)[(1,0),\cdot]+\de^2\varphi_s(x)[y,\cdot]\big].
\end{multline*}
Moreover, as $\partial_yG(x,y)$ and $\partial_x G(x,y)$ are homogeneous in $y$ of degree $1$ and $2$ respectively, using last equation, recalling that $\{\xi_n\}$ is bounded in the $C^0$
norm and the fact that all the involved operators are uniformly bounded in norm
 we get
\bal
\lefteqn{\left|\inte\partial_x\tilde{G}_s(\xi_n,\dot  \xi_n)[\xi_n-\xi_m]\de s\right|\leq}&\nonumber\\
&K_6\inte (1+|\dot\xi_n|^2)|\xi_n-\xi_m|\de s+K_7\inte (1+|\dot\xi_n|^2)^{1/2}(1+|\dot\xi_n|)|\xi_n-\xi_m|\de s,
\label{tispiezzoindue}
\eal
for some positive constants $K_6$ and $K_7$. As by \eqref{xil2bound}, $\{\dot\xi_n\}$ is  bounded in the $L^2$ norm,   and $\{\xi_n\}$ is Cauchy in the $C^0$ norm, the right-hand side in \eqref{tispiezzoindue} and therefore  the first integral in \eqref{carta} goes to $0$ as $m,n\to \infty$.

Now we change the role of $\xi_n$ and $\xi_m$ considering $\de \tilde{J}
(\xi_m)[\xi_n-\xi_m]$.
Proceeding as in \eqref{byps}, we see that $\de \tilde{J}(\xi_m)[\xi_n-\xi_m]\to 0$, as
$m,n\to\infty$. Therefore
\beq\label{changer}
\inte\big(\partial_y\tilde{G}_s(\xi_n,\dot \xi_n)[\dot \xi_n-\dot \xi_m]
-\partial_y\tilde{G}_s(\xi_m,\dot \xi_m)[\dot \xi_n-\dot \xi_m]\big)\de s\To 0, \ \
\text{as $m,n\to\infty.$}\eeq
Since $\partial_y\tilde G_s(x,y)[\cdot]=\partial_yG\big(\varphi(s,x),\de\varphi(s,x)[(1,y)]\big)\big[\de\varphi_s(x)[\cdot]\big]$ and $\xi_n\to\bar\xi$ uniformly, using the facts that $\partial_yG(x,y)$ is continuous on $TM$ and positively homogeneous  of degree $1$ in $y$, that a continuous function on a compact set is uniformly continuous and that $\{\dot\xi_n\}$ is uniformly bounded in the $L^2$ norm, the limit  \eqref{changer} gives also
\beq\label{richanger}
\inte\big(\partial_y\tilde{G}_s(\bar \xi,\dot \xi_n)-\partial_y\tilde{G}_s(\bar \xi,\dot
\xi_m)\big)[\dot \xi_n-\dot \xi_m]\de s
\To 0, \ \ \text{as $m,n\to\infty$.}\eeq
Let us define $\delta_i(s)=\de \varphi(s,\bar{\xi}(s))[(1,\dot\xi_i(s))]$ for $i\in\N$ and $s\in[0,1]$ and  the following subsets of the interval $[0,1]$.  Let $A_i\subset [0,1]$ be the support of the $L^2$ function  $[0,1]\ni s\to|\delta_i(s)|=\big(h(\bar{\xi}(s))[\delta_i(s),\delta_i(s)]\big)^{1/2}\in\R$ for $i\in\N$ and choose
\baln
&B_{nm}=\left\{t\in A_n\cap A_m\ \left|\ \frac{\delta_m}{|\delta_m|}=- \frac{\delta_n}{|\delta_n|}\text{ a.e.} \right.\right\},\\
&C_{nm}=(A_n\cup A_m)\setminus B_{nm},\\
&D_{nm}=[0,1]\setminus(A_n\cup A_m).
\ealn
Moreover, we assume that $\frac{\delta_m}{|\delta_m|}\not=- \frac{\delta_n}{|\delta_n|}\text{ a.e.}$ in $A_n\cap A_m\setminus B_{nm}$.
We observe that the interval $[0,1]$ is the union of the sets $B_{nm}$, $C_{nm}$ and
$D_{nm}$, for every $n$ and $m$; moreover on $B_{nm}$ we have $\delta_m=-\lambda_{nm}\delta_n$, with $\lambda_{nm}=\frac{|\delta_m|}{|\delta_n|}$. The subsets $B_{nm}$ and $D_{nm}$ are precisely the instants where the mean value theorem cannot be applied because of the lack of smoothness of $G$ on the null section.
Applying the mean value theorem for every $s\in C_{nm}$ and  using the fact that
\[\partial_{yy}\tilde{G}_s(x,y)[\cdot,\cdot]=\partial_{yy}G\big(\varphi(s,x),\de\varphi(s,x)[(1,y)]\big)\big[\de\varphi_s(x)[\cdot],\de\varphi_s(x)[\cdot]\big]
\]
is positive definite and $\partial_{yy} G(x,y)$ is positive
homogeneous of degree $0$ in $y$, we get the existence of a positive constant
$K_8$ such that
\bml
\int_{C_{nm}} \big(\partial_y\tilde{G}_s(\bar \xi,\dot \xi_n)-\partial_y\tilde{G}_s(\bar \xi,\dot
\xi_m)\big)[\dot \xi_n-\dot \xi_m]\de s\\
=\int_{C_{nm}}\partial_{yy}\tilde{G}_s\big(\bar \xi,\vartheta\dot \xi_n+(1-\vartheta)\dot
\xi_m\big)[\dot \xi_n-\dot \xi_m,\dot \xi_n-\dot \xi_m]\de s\\\geq K_8\int_{C_{nm}}|\dot \xi_n-\dot \xi_m|^2\de s.\label{meanvalue}
\eml
where $\vartheta\colon C_{nm}\to\R$ is a function assuming values in
$[0,1]$.
We pass now to estimate the functions
$\big(\partial_y\tilde{G}_s(\bar{\xi},\dot\xi_n)-\partial_y\tilde{G}_s(\bar{\xi},\dot\xi_m)\big)[\dot\xi_n-\dot\xi_m]$
over the subsets $B_{nm}$.
To this end, we observe that
\baln
\lefteqn{\big(\partial_y\tilde{G}_s(\bar\xi,\dot\xi_n)-\partial_y\tilde{G}_s(\bar\xi,\dot\xi_m)\big)[\dot\xi_n-\dot\xi_m]}&\\
&=\Big(\partial_yG\big(\varphi_s(\bar\xi),\de\varphi(s,\bar\xi)[(1,\dot\xi_n)]\big)\\
&\quad-\partial_yG\big(\varphi_s(\bar\xi),\de\varphi(s,\bar\xi)[(1,\dot\xi_m)]\big)\Big)\big[\de\varphi_s(\bar\xi)[\dot\xi_n-\dot\xi_m]\big]\\
&=\Big(\partial_yG\big(\varphi_s(\bar\xi),\de\varphi(s,\bar\xi)[(1,\dot\xi_n)]\big)\\
&\quad-\partial_yG(\varphi_s(\bar\xi),\de\varphi(s,\bar\xi)[(1,\dot\xi_m)]\big)\Big)\big[\de\varphi(s,\bar\xi)[(0,\dot\xi_n)]-\de\varphi(s,\bar\xi)[(0,\dot\xi_m)]\big]\\
&=\Big(\partial_yG\big (\varphi_s(\bar\xi),\de\varphi(s,\bar\xi)[(1,\dot\xi_n)]\big)\\
&\quad-\partial_yG\big (\varphi_s(\bar\xi),\de\varphi(s,\bar\xi)[(1,\dot\xi_m)]\big)\Big)\big[\de\varphi(s,\bar\xi)[(1,\dot\xi_n)]-\de\varphi(s,\bar\xi)[(1,\dot\xi_m)]\big].
\ealn
Therefore, recalling that  $\delta_m=-\lambda_{nm}\delta_{ n}$ over the subsets $B_{nm}$, we get
\baln
&\int_{B_{nm}}
\big(\partial_y\tilde{G}_s(\bar{\xi},\dot\xi_n)-\partial_y\tilde{G}_s(
\bar{\xi},\dot\xi_m)\big)[\dot\xi_n-\dot\xi_m]\de s\\
&=\int_{B_{nm}} \big(\partial_yG(\varphi_s(\bar{\xi}),\delta_n)-\partial_yG(
\varphi_s(\bar{\xi}),\delta_m)\big)[\delta_n-\delta_m]\de s\\
&=\int_{B_{nm}}(1+\lambda_{nm}) \partial_yG(\varphi_s(\bar{\xi}),\delta_n)[\delta_n]\de
s\\
&\quad+\int_{B_{nm}}\left(1+\frac{1}{\lambda_{nm}}\right)\partial_yG(\varphi_s(\bar{\xi}),\delta_m)[\delta_m]\de s.
\ealn
By Euler's theorem the above integrals are equal to
\[
\int_{B_{nm}}2(1+\lambda_{nm}) G(\varphi_s(\bar{\xi}),\delta_n)\de s  +\int_{B_{nm}}2\left(1+\frac{1}
{\lambda_{nm}}\right)G(\varphi_s(\bar{\xi}),\delta_m)\de s
\]
and then, by homogeneity, we get
\beq\label{onBnm}
\begin{split}
\int_{B_{nm}}2 (1+\lambda_{nm})G(&\varphi_s(\bar \xi),\delta_n)\de s  +\int_{B_{nm}}2\left(1+\frac{1}
{\lambda_{nm}}\right)G(\varphi_s(\bar \xi),\delta_m)\de s
\\
&\geq K_9\left(\int_{B_{nm}}|\delta_n|^2\de s+\int_{B_{nm}}|\delta_m|^2\de s\right),\\
\end{split}
\eeq
 where $K_9$ is a positive constant independent of $n$ and $m$.
Now observe that, by linearity
\[|\delta_n-\delta_m|^2=|\de\varphi(s,\bar\xi)[(0,\dot\xi_n-\dot\xi_m)]|^2\geq K_{10}|\dot\xi_n-\dot\xi_m|^2,\]
where $K_{10}$ is  the minimum value of the function $|\de \varphi(s,\bar\xi(s))[(0,v)]|^2$ over the compact set $[0,1]\times S^{r-1}$, where $S^{r-1}$ is the $(r-1)$-dimensional sphere. Therefore we obtain
\bml\label{fiu}\int_{B_{nm}}\!\!|\dot \xi_n-\dot\xi_m|^2\de s\leq \frac{1}{K_{10}}\int_{B_{nm}}\!\!|\delta_n-\delta_m|^2\de s\\\leq\frac{2}{K_{10}}\left(\int_{B_{nm}}|\delta_n|^2\de s+\int_{B_{nm}}|\delta_m|^2\de s\right).\eml
Over the subset $D_{nm}$ both $\delta_n$ and $\delta_m$ are zero, hence
\beq\label{fiuu}
\int_{D_{nm}}|\dot \xi_n-\dot \xi_m|^2\de
s\leq \frac{1}{K_{10}}\int_{D_{nm}}|\delta_n-\delta_m|^2\de
s=0,
\eeq
for all $n$ and $m$.
From \eqref{richanger}, summing up \eqref{meanvalue}, \eqref{onBnm}, \eqref{fiu} and \eqref{fiuu} and recalling that the interval $[0,1]=B_{nm}\cup C_{nm}\cup D_{nm}$, we finally get
\[
\inte|\dot \xi_n-\dot \xi_m|^2\de s\To 0\]
as $n,m\to\infty$.
\qed\end{proof}
With the Palais-Smale condition  in hand,  infinite dimensional
Lusternik and Schnirelman theory becomes available (see
\cite{Palais66}); so we can  obtain  existence and  multiplicity
results about the number of critical points of $J$, depending on $N$
and the topology of $M$, for example in the non-contractible case.   We consider
here the case of geodesics joining two different submanifolds of
$M$ (compare also with \cite[Theorem 6]{KoKrVa04}).
 \bpr\label{multiplicity} 
Let $(M,F)$ be a forward or 
backward complete Finsler manifold and let $M_1$ and $M_2$ be two
closed submanifolds of $M$ such that $M_1$ or $M_2$ is compact. 
Then there exists a geodesic $\gamma$
connecting $M_1$ and $M_2$ and satisfying \eqref{orto}. Moreover, if
the manifold $M$ is non-contractible and $M_1$, $M_2$ are contractible 
then there exist
infinitely many geodesics $\gamma_n$ connecting $M_1$ and $M_2$,
satisfying \eqref{orto} and such that
$\lim_nJ(\gamma_n)=+\infty$ (according to Theorem~\ref{ps}, in the forward case such geodesics start from the compact submanifold while, in the backward case, they arrive to it). 
\epr
\begin{proof}
Existence is a standard application of the Deformation Lemma (see \cite{Palais66}). For the multiplicity result we recall that, given a topological space $X$, the
Lusternik-Schnirelman category of $X$,  is a homotopy invariant
defined as the minimum number, denoted by $\cat X$, of closed
contractible subsets of $X$ which cover $X$. Let $C^0_{M_1\times
M_2}(M)$ be the space of the continuous curves having endpoints in
$M_1\times M_2$. The inclusion of $\Lambda_{M_1\times M_2}(M)$ in
$C^0_{M_1\times M_2}(M)$ is a homotopy equivalence (see
\cite[Theorem 1.3]{Grove73}). Let $\Omega(M)$ be the space of based
loops in $M$. Since $M_1$ and $M_2$ are contractible,
$C^0_{M_1\times M_2}(M)$ has the same homotopy type  as
$M_1\times M_2\times \Omega(M)$, moreover $\cat \Omega(M)=\infty$ (see
\cite[Proposition 3.2 and Corollary 3.2]{FadHus91}), hence
\[
\mathrm{cat}(\Lambda_{M_1\times M_2}(M))=\infty
\]
as well. By Theorem 7.2 of \cite{Palais66} and Theorem~\ref{ps},
$J$ has infinitely many critical points $\gamma_n$ which are
geodesics connecting $M_1$ to $M_2$ and satisfying  \eqref{orto}.
Finally, $\sup_{n\in\N}J(\gamma_n)=+\infty$ otherwise would be
possible to retract  the manifold $\Lambda_{M_1\times M_2}(M)$ onto
a sublevel of the functional $J$. This would be a contradiction,
since the sublevels of a $C^1$ functional  defined on a Banach
manifold, bounded from below and satisfying the Palais-Smale
condition have finite Lusternik-Schnirelman category.
\qed\end{proof}

\bere We point out that the above multiplicity result does not
guarantee, in general, that the infinitely many geodesics are
geometrically distinct (they might cover the
same closed geodesic, as on the standard  sphere). \ere
\bere For the two endpoints boundary conditions, the
above multiplicity result can also be obtained by using Morse theory
and a finite dimensional approximation of the path space
$\Lambda_{\{p\}\times\{q\}}(M)$, see \cite[Ch. III \S\S
16,17]{Milnor63} (cf. \cite{Shen01} for the Finsler case). However,
for general boundary conditions the infinite dimensional approach is very useful.
In particular for periodic boundary conditions, in contrast to the finite dimensional
approximation, the free loop space carries, for a non-reversible
Finsler metric, a canonical $S^1$-action leaving the energy
functional invariant. \ere
\section{The Fermat metric}\label{appli}
In this section we present some applications of Finsler Geometry to
the study of the causal structure of a conformally stationary
spacetime. We first recall the definition of a Finsler manifold of
Randers type, then we introduce a  Randers metric, that we call {\em
Fermat metric}, which is related to the Fermat principle for
lightlike geodesics in a conformally stationary spacetime.
\subsection{Randers metrics} Let $h$ be a Riemannian tensor
and $\omega$ be a one-form on $M$.  A Randers metric $F$ is defined
as follows: \beq\label{randers}
F(x,y)=\sqrt{h(x)[y,y]}+\omega(x)[y],\quad\quad\quad
\|\omega\|_x<1,\eeq where  $\|\omega\|_x= \sup_{v\in T_x M\setminus
0}\frac{|\omega(x)[v]|}{\sqrt{h(x)[v,v]}}$. Remarkably enough, the
condition $\|\omega\|_x<1$ for all $x\in M$, not only implies that
$F$ is positive but also that it has fiberwise  strongly convex
square (see  \cite[\S 11.1]{BaChSh00}).

\begin{remark}\label{completo}
We observe that if the Riemannian metric $(M,h)$ is complete and \beq\label{superscemo}
\|\omega\|\colon=\sup_{x\in M}\|\omega\|_x<1,
\eeq
the Randers manifold $(M,F)$ is forward  and  backward  complete.
In fact, let $\{x_n\}$ be, for instance, a forward  Cauchy sequence for $(M,F)$,
then for any $\eps>0$ there exists $\nu\in\N$ such that for all $i,\ j\in\N$ with $\nu\leq i\leq j$, $\dist(x_i,x_j)<\eps$.
By definition of distance, there exists a curve $\gamma_{ij}$ connecting $x_i$ to $x_j$, such that
\[\eps\!>\!\int_{\gamma_{ij}}\!\!\! F(\gamma_{ij},\dot\gamma_{ij})\geq (1-\|\omega\|)\!\!\int_{\gamma_{ij}}\!\! \sqrt{h(\gamma_{ij})[\dot\gamma_{ij},\dot\gamma_{ij}]}\geq (1-\|\omega\|)\distg(x_i,x_j),\]
where \distg is the distance associated to the Riemannian metric $h$. Being $(M,h)$ complete, $\{x_n\}$ converges.
\end{remark}

\subsection{The Fermat metric of a conformally standard stationary spacetime}\label{lightlike}
A Fermat principle in General Relativity is  a variational
characterization of the light rays  joining an event with the
worldline of an observer in the spacetime. A spacetime is given by a
Lorentzian manifold whose metric tensor satisfies the Einstein
equations together with a time orientation, while  light rays are given by the lightlike geodesics of
the Lorentzian manifold. In the recent years there has  been a great
amount of work  about the Fermat principle in General Relativity,
because it allows one to obtain a mathematical description of the {\it
gravitational lens effect} in Astrophysics, see
\cite{GiMaPi02,Perlic04}.

A Lorentzian manifold $({\mathcal M},g)$ is a smooth, 
connected spacelike manifold ${\mathcal M}$ endowed with a symmetric
non-degenerate tensor field $g$ of type $(0,2)$ having index $1$. A
geodesic of $({\mathcal M},g)$ is a smooth curve
$\gamma\colon[a,b]\to \mathcal M$ satisfying the equation
$\nabla_{\gamma}\dot\gamma=0$, where $\nabla_{\gamma}$ is the
covariant derivative along $\gamma$ associated to the Levi-Civita
connection of the metric $g$ (we refer to \cite{BeErEa96} for all
the needed background material on Lorentzian geometry). As in the
Riemannian case, a geodesic has to satisfy the conservation law
$g(\gamma)[\dot\gamma,\dot\gamma]=E_{\gamma}=\mathrm{const.}$,  which
corresponds to energy conservation in  Lagrangian mechanics.
According to the sign of $E_{\gamma}$, a  geodesic is said {\em
timelike} if $E_{\gamma}<0$, {\em lightlike} if $E_{\gamma}=0$, {\em
spacelike} if $E_{\gamma}>0$ or $\dot\gamma(s)=0$ for all
$s\in[a,b]$. This  partition  of the set of  geodesics is  known  as
the {\em the causal character} of a geodesic. Such a terminology is
used also for any vector in any tangent space and for any piecewise
smooth curve if and only if its  tangent vector field has the same
character at any point where it is defined.

A time  orientation on a Lorentzian manifold is determined by a
timelike vector field $Y$ on ${\mathcal M}$, i.e. for any $p\in
{\mathcal M}$, $Y(p)$ is a timelike vector. A piecewise smooth
{lightlike)  curve $\gamma: [a,b]\to \mathcal M$ is said to be {\em
future-pointing} (resp. {\em past-pointing}) if $g(\gamma(s))[\dot
\gamma(s),Y(\gamma(s))]<0$ (resp. $g(\gamma(s))[\dot
\gamma(s),Y(\gamma(s))]>0$) for all $s\in[a,b]$ where
$\dot\gamma(s)$ is defined. The notion of being future-pointing
(resp. past-pointing) and non-spacelike can be extended to a
continuous curve $\gamma\colon[a,b]\to {\mathcal M}$ requiring that
for any $s_0\in[a,b]$ there is a convex normal neighborhood
$U\subset {\mathcal M}$ of $\gamma(s_0)$ and  an interval
$J\subset[a,b]$ containing $s_0$ such that for any $s_1,\ s_2\in J$,
with $s_1<s_2$, a smooth future-pointing (resp. past-pointing) curve
connecting $\gamma(s_1)$ to $\gamma(s_2)$ and contained in $U$
exists. \label{causal} From now on, non-spacelike curves are assumed
to be future-pointing.

A {\em conformally standard stationary Lorentzian} manifold is a  manifold
${\mathcal M}$ which splits as a product ${\mathcal M}={\mathcal
M}_0\times \R$, where ${\mathcal M}_0$ is endowed with a Riemannian
metric $g_0$,  with a vector field $\delta$ and a 
positive function $\beta$. Moreover, there exists a positive function  $\varphi$ on ${\mathcal M}$, such that the
Lorentzian metric $g$ on ${\mathcal M}$ is given by
\begin{equation}\label{l2}
g(x,t)[(y,\tau),(y,\tau)]=\varphi(x,t)\big(g_0(x)[y,y]+2g_0(x)[\delta(x),y]\tau
-\beta(x)\tau^2\big),
\end{equation}
for any $(x,t)\in {\mathcal M}_0\times \R$ and  $(y,\tau)\in
T_x{\mathcal M}_0\times\R$.   A conformally standard stationary Lorentzian manifold is time oriented by the
timelike Killing vector field $\partial_t$ and a piecewise smooth
non-spacelike curve $\gamma=(x,t)$ is future-pointing iff $\dot t>0$
where $\dot\gamma$ exists.

Since lightlike geodesics and causal properties  - as global hyperbolicity - of a conformally stationary spacetime are invariant  under conformal changes of the metric tensor $g$ (see for example \cite{BeErEa96,Piccio97}), we can assume that $g$ is given by $g/\varphi(x,t)$. Indeed, rather than the metric in \eqref{l2}, we will consider the metric
\begin{equation}\label{l}
g(x,t)[(y,\tau),(y,\tau)]=g_0(x)[y,y]+2g_0(x)[\delta(x),y]\tau
-\beta(x)\tau^2.
\end{equation}

We introduce now the Fermat metric associated to a standard
stationary Lorentzian manifold. Let $z_0 = (x_0,t_0) \in {\mathcal M}$ be an event in ${\mathcal M}$ and let
$\R\ni s\rightarrow\gamma(s) = (x_1,s)\in{\mathcal M}$
be a timelike vertical curve, that is an integral curve of the timelike vector field $\partial_t$.
Let $[0,1]\ni s\rightarrow z(s) = (x(s),t(s))\in{\mathcal M}$ be a  lightlike curve joining $z_0$ and $\gamma$.
Concretely the lightlike curve $z$ satisfies
\begin{equation}\label{light}
g_0(x)[\dot x,\dot x]+2g_0(x)[\delta(x),\dot x]\dot t-\beta(x)\dot t^2=0,
\end{equation}
and the boundary conditions $x(0) = x_0$, $x(1) = x_1$, $t(0)=t_0$.
The {\it arrival time} $T(z)$ of the lightlike curve $z$ is given by $t(1)$. From \eqref{light},  assuming that the lightlike curve is future-pointing,
solving with respect to $\dot t$ and integrating we obtain:
\begin{equation}\label{tempo}
t(s)=t_0+\!\int_0^s\big(\tig(x)[\delta(x),\dot x]+
\sqrt{\tig(x)[\delta(x),\dot x]^2+\tig(x)[\dot x,\dot x]}\big)\de v,
\end{equation}
where $\tig=g_0/\beta$.
So the arrival time $T(z)$ is given by
\beq\label{Arrivaltime}
T(z)=t_0+\!\inte \big(\tig(x)[\delta(x),\dot x]+
\sqrt{\tig(x)[\delta(x),\dot x]^2+\tig(x)[\dot x,\dot x]}\big)\de s.
\eeq
\bd\label{Fermatmetric}
The {\em Fermat metric} associated to a standard stationary Lorentzian manifold $(\mathcal M,g)$ ($g$ as in \eqref{l}) is the Randers metric $F$ on $\mathcal M_0$  given by
\[
F(x,y)=\tig(x)[\delta(x),y]+
\sqrt{\tig(x)[\delta(x),y]^2+\tig(x)[y,y]}\]
for every $(x,y)\in T\mathcal M_0$, being $\tig=g_0/\beta$ (cf. \eqref{randers}; here the Riemannian metric $h$ is given by $h(x)[y,y]=\tig(x)[\delta(x),y]^2+\tig(x)[y,y]$).
\ed
\bere
The Fermat metric associated to a conformally standard stationary spacetime as in \eqref{l2} will be the Fermat metric associated to the standard stationary spacetime $(\mathcal M,g/\varphi(x,t))$.
\ere
\bere
Observe that the arrival time $T(z)$ of a future-pointing lightlike curve $z = (x,t)$ is equal, up to  the initial instant of time $t_0$, to the length of its spatial projection $x$ with respect to the Fermat metric $F$.
\ere
The Fermat metric $F$ allows one to reduce Fermat's principle for light rays on a  standard stationary spacetime to a variational principle  involving only the spatial projections of  the lightlike curves.
We recall that the relativistic Fermat principle for
lightlike geodesics states that among all lightlike curves $z\colon[0,1]\to \mathcal M$ connecting some
event $p\in \mathcal M$ with some timelike curve $\gamma$ on $\mathcal M$, lightlike geodesics
are, up to reparameterization, critical points of the arrival time,
which is the functional $z\mapsto \gamma^{-1}(z(1))$.
The property of lightlike geodesics (light
rays) of  being stationary points of the arrival time is classically
known as {\em Fermat's principle}. The first one to formulate
Fermat's principle in General Relativity in the above generality
was I. Kovner in \cite{Kovner90}, but a rigorous proof was
given by V. Perlick in \cite{Perlic90}. Some special versions of
Fermat's principle for static, stationary and conformally stationary
Lorentzian manifolds, contained in several books about General Relativity, can be
deduced from the above general version.
The Finslerian reduction of the principle for a standard stationary spacetime consists
in proving that  a future-pointing lightlike curve $[0,1]\ni s\to(\tilde x(s),\tilde t(s))\in {\mathcal M}$ joining $(x_0,t_0)$ with  $\gamma(s)=(x_1,s)$ is a lightlike geodesic of $({\mathcal M},g/\beta)$ (and up to reparametrization of $({\mathcal M},g)$) if and only if its
spatial component $\tilde x$ is a geodesic of the Fermat metric.
\bt[Fermat's principle]
Let $(\mathcal M,g)$ be a standard stationary spacetime and $(x_0,t_0)\in\mathcal M$,  $\R\ni s\to\gamma(s)=(x_1,s)\in \mathcal M$, $x_1\in\mathcal M_0$.  A curve $[0,1]\ni x\to z(s)=(x(s),t(s))\in\mathcal M$ is a future-pointing lightlike geodesic of $({\mathcal M},g/\beta)$ if and only if $x(s)$ is a geodesic for the Fermat metric $F$, parameterized to have constant Riemannian speed $h(x)[\dot x,\dot x]=\tig(x)[\delta(x),\dot x]^2+\tig(x)[\dot x,\dot x]$, and $t(s)$ is given by \eqref{tempo}.
\et
\begin{proof}
Using the Levi-Civita connection $\nabla$ of the metric  $\tig$, the Euler-Lagrange equations  of  the functional \eqref{Arrivaltime}
can be written as
\begin{equation}\label{eulerlagrange}-\nabla_{\dot x}\left(\frac{\dot x+\tig[\delta,\dot x]\delta}{\sqrt{h[\dot x,\dot x]}}\right)+\frac{\tig[\delta,\dot x] (\nabla\delta)^*[\dot x]}{\sqrt{h[\dot x,\dot x]}}+ (\nabla\delta)^*[\dot x]- \nabla\delta[\dot x]=0,
  \end{equation}
where $(\nabla\delta)^*$ is the
adjoint with respect to $\tig$ of $\nabla\delta$ and $(\nabla
\delta)[\dot x]=\nabla_{\dot x}\delta$.
Hence if $x$ is parameterized to have  constant Riemannian speed $h[\dot x,\dot x]$,  we get:
\bal
\nabla_{\dot x}\dot  x&=-\nabla_{\dot x}\left(\tig[\delta,\dot x]\delta\right)+\tig[\delta,\dot x] (\nabla\delta)^*[\dot x]+\sqrt{h[\dot x,\dot x]}\left((\nabla\delta)^*[\dot x]- \nabla\delta[\dot x]\right)\nonumber\\
&=-\frac{\de}{\de s}\left(\tig[\delta,\dot x]\right)\delta+\tig[\delta,\dot x]\left( (\nabla\delta)^*[\dot x]- \nabla\delta[\dot x]\right)\nonumber\\
&\quad\quad+\sqrt{h[\dot x,\dot x]}\left((\nabla\delta)^*[\dot x]- \nabla\delta[\dot x]\right)\nonumber\\
&=F( x,\dot x)\Omega[\dot x]-\frac{\de}{\de s}\left(\tig[\delta,\dot
x]\right)\delta,\label{hconst}\eal where $\Omega[\dot
x]=(\nabla\delta)^*[\dot x]-( \nabla\delta)[\dot x]$.  Lightlike geodesics of $(\mathcal M,g/\beta)$ are  critical points of the
energy functional 
\[(x,t)\mapsto \frac1 2\int_0^1(\tig[\dot x,\dot
x]+2\tig[\delta,\dot x]\dot t-\dot t^2)\de s,\] so they satisfy the
Euler-Lagrange equations
\beq\label{EL}
 \begin{cases}
\nabla_{\dot x}\dot x =\dot t\Omega[\dot x]-\frac{\de\dot t}{\de s}\delta,\\
\tig[\delta,\dot x]+C=\dot t,
 \end{cases}
\eeq
where $C$ is a constant. By the second equation in \eqref{EL},   $\frac{\de \dot t}{\de s}=\frac{\de}{\de s}(\tig[\delta,\dot x])$ and recalling that a future-pointing lightlike curve has to satisfy
the   equation
\beq\label{fplc}
\dot t=F(x,\dot x),
\eeq
we get \eqref{hconst}. Finally integrating \eqref{fplc} we get that $t(s)$ is given by \eqref{tempo}. The reciprocal is analogous.
\qed\end{proof}
\bere
We point out that the name Fermat metric has been used in some paper to denote the Riemannian metric $\tilde g_0$ (see \cite[\S 4.2]{Perlic04} and the references therein). We think that
our definition is more appropriate because, as for the Fermat principle in classical optics, arrival times of
lightlike curves and in particular light rays are measured as lengths with respect to a metric, in this case a Finsler one.
\ere
We shall see now how the Fermat metric has not only a clear variational meaning, but it plays a basic role also in the study of causal properties of a conformally standard stationary spacetime.
We recall some basic  definitions and  properties about causality  (our main references about that are  \cite{BeErEa96,HawEll73}).
A Lorentzian manifold $({\mathcal M},g)$ is said {\em strongly causal} if every $p\in {\mathcal M}$ has arbitrarily small neighborhoods
such that no non-spacelike curve that leaves one of these neighborhoods ever returns.
A non-spacelike curve $\gamma\colon(a,b)\to {\mathcal M}$ is said {\em future inextendible} (resp. {\em past
 inextendible})
if the limit $\lim_{s\to b^-}\gamma(s)$ (resp. $\lim_{s\to a^+}\gamma(s)$) does not exist.
It is said {\em inextendible} if it is both future and past inextendible.
Two non-spacelike continuous curves are considered equivalent if one is the reparameterization of the other.
Henceforth, whenever the domain of the parameter is not specified, we will be regarding the equivalence class of the curve.
For any $p\in {\mathcal M}$, let $J^+(p)\subset {\mathcal M}$ (resp. $J^-(p)\subset {\mathcal M}$) be the subset of the points $q$ in ${\mathcal M}$
such that there exists a non-spacelike curve $\gamma\colon[a,b]\to {\mathcal M}$ with $\gamma(a)=p$ and $\gamma(b)=q$
(resp. $\gamma(a)=q$ and $\gamma(b)=p$).
A Lorentzian manifold $({\mathcal M},g)$ is said   {\em globally hyperbolic} if it
admits a {\em Cauchy surface} i.e. a subset $S$ which every inextendible  timelike curve intersects exactly once.
It can be proved that $({\mathcal M},g)$ is globally hyperbolic if and only if it is strongly causal and for all $p,\ q\in {\mathcal M}$ the set
$J^+(p)\cap J^-(q)$ is compact (see \cite[Proposition 6.6.3 and Proposition 6.6.8]{HawEll73}).

For future references, we show here, in the case of a conformally standard stationary Lorentzian manifold, a fact that is cited,
without proof, in several references (see for example \cite[p. 213]{HawEll73}).
\begin{lemma}\label{ligthwise}
Let $({\mathcal M},g)$ be a conformally standard stationary Lorentzian manifold and $p,q\in \mathcal M$ two  causally connected points.
Then there is  a piecewise lightlike geodesic connecting $p$ and $q$.
\end{lemma}
\begin{proof}
As causality is invariant by conformal transformations we can assume that the metric is standard stationary as in \eqref{l}.
Let $\gamma:[0,1]\to {\mathcal M}$ be a non-spacelike curve joining $p$ and $q$ given by $\gamma(s)=(x(s),t(s))$.
To find a piecewise smooth lightlike geodesic connecting $p$ and $q$ is equivalent to finding a piecewise
smooth Fermat geodesic joining $x(0)$ and $x(1)$  and having length equal to $t(1)-t(0)$.
As convex neighborhoods always exist in Finsler geometry (see
\cite{Whiteh32}), the support of $x$ can be covered by a finite number of them. So we can assume, without loss of generality,
 that $x(0)$ and $x(1)$ are in the same convex neighborhood.
 We will show that there exist piecewise smooth geodesics from $x(0)$ to $x(1)$ having length
 $s$ for every $s\geq\dist(x(0),x(1))$, and then the result follows, since $ t(1)-t(0)=\inte\dot t(s)\de s\geq\inte F(x,\dot x)\de s\geq \dist(x(0),x(1))$,
 where the first inequality comes from the inequality
 $g_0(x)[\dot x,\dot x]+2g_0(x)[\delta(x),\dot x]\dot t -\beta(x)\dot t^2\leq 0$ which says that
 $\gamma$ is a non-spacelike curve.
 First,  observe that there is a minimal geodesic from $x(0)$ to $x(1)$ with Fermat length equal to $\dist(x(0),x(1))$,
 because they are contained in a convex neighborhood.
 Then we can choose two sequences  of points $\{x_i\}$ and $\{y_j\}$ in such a way that
 the distance between one element of the first sequence and another of the second one is always bigger than a small enough $\varepsilon>0$. 
 Making a sufficient number of ``zig zags'', the length of the piecewise geodesic can  be made as big as needed.
 In the first ``zig zag'' where the piecewise geodesic length becomes  bigger than $s$,
 we can move  back the last point along the last piece of geodesic.
 As the variation of the length is continuous, we can construct, in  this way, a piecewise geodesic with length $s$.
\qed\end{proof}
In the next theorem, we show that on a conformally standard stationary Lorentzian manifold,
global hyperbolicity is  strictly related to the Fermat metric completeness.
To the authors' knowledge,   this link between global hyperbolicity and the completeness of the Fermat metric   does not appear elsewhere
in literature. We are going to use the following notation for $p_0=(x_0,t_0)\in {\mathcal M}$:
${\rm C}^+(p_0,\mu)=\cup_{s\in[0,\mu)} \bar{B}^+_s(x_0)\times\{t_0+s\}$ and ${\rm C}^-(p_0,\mu)=\cup_{s\in[0,\mu)} \bar{B}^-_s(x_0)\times\{t_0-s\}$,
where $\bar{B}^{\pm}_s(x_0)$ is the closure of $B^{\pm}_s(x_0)$ in ${\mathcal M}_0$.

\begin{theorem}\label{gh}
Let $({\mathcal M},g)$ be a conformally standard stationary Lorentzian manifold and let $\bar{t}\in\R$. Then the following propositions hold:
\begin{itemize}
\item[(1)]  if the Fermat metric on ${\mathcal M}_0$ defined in \ref{Fermatmetric}
is forward (or resp. backward) complete then
$J^+(p_0)={\rm C}^+(p_0,+\infty)$ and $J^-(p_0)={\rm C}^-(p_0,+\infty)$ for every $p_0=(x_0,t_0)\in L$, the balls $\bar{B}^+_s(x_0)$
(resp. $\bar{B}^-_s(x_0)$) are compact and $({\mathcal M},g)$ is globally hyperbolic;
\item[(2)]  if $({\mathcal M},g)$ is globally hyperbolic with Cauchy surface $S={\mathcal M}_0\times \{\bar{t}\}$ then the Fermat metric
on ${\mathcal M}_0$ is forward and backward complete.
\end{itemize}
\end{theorem}

\begin{proof}
Again we can assume that $g$ is as in \eqref{l}.
 We begin with proving $(1)$, assuming that $F$ is forward complete (the proof in the backward case is analogous).
 Compactness of the balls $\bar{B}^+_s(x_0)$ is a consequence of the Finslerian Hopf-Rinow theorem.
 Now assume that $(x_1,t_1)\in \bar{B}^+_s(x_0)\times\{t_0+s\}$ for a certain $s\in [t_0,+\infty)$.
 By applying  the Finslerian Hopf-Rinow theorem  we can choose a Finslerian minimal geodesic $x$ from $x_0$ to $x_1$ with speed equal to $1$ and
 length not greater  than $s$.
 Considering the lightlike geodesic $\lambda\to(x(\lambda), \lambda+t_0)$ with $\lambda\in  [0,L(x)] $, where $L(x)$ is the Fermat length of $x$,
 and then the timelike curve $\lambda\to  (x(L(x)),\lambda)$ with $\lambda\in[ t_0+L(x),t_0+s]$,
 we see that $(x_1,t_1)\in J^+(p_0)$.
 If $q=(x_1,t_1)\in J^+(p_0)$, then by Lemma~\ref{ligthwise} there exists a piecewise smooth lightlike geodesic
 $\gamma(s)=(x(s),t(s))$ which connects $p_0$ to $q$, such that $\dist(x_0,x_1)\leq  L(x)$, hence $q\in \bar B^+_{ L(x)}(x_0)\times\{t_0+L(x)\}$.
Analogously one can prove the other equality.
 Now since   $({\mathcal M},g)$ admits the coordinate $t$ as a global time function  it is stably causal and then strongly causal
 (see for instance \cite[p. 64 and p. 73]{BeErEa96}).
 Furthermore, if $p=(\bar{x},\bar{t})$ and $q=(\tilde{x},\tilde{t})$ are points in ${\mathcal M}$, then we can assume that $\tilde t>\bar t$
 otherwise $J^+(p)\cap J^-(q)$ is empty. Moreover we have
\[J^+(p)\cap J^-(q)=\bigcup_{s\in[0,1]}\left(\bar B^+_{s(\tilde t-\bar t)}(\bar x)\cap
\bar B^-_{(1-s)(\tilde t-\bar t)}(\tilde x)\right)\times\{\bar t+s(\tilde t - \bar t)\},
\]
which is compact or empty.
This can be shown as follows. Take a sequence $\{(x_n,t_n)\}\subset J^+(p)\cap J^-(q)$; as $\{t_n\}$ moves in a compact set,
we can extract a convergent subsequence.
If $\bar{r}=\sup_n \{t_n -\bar{t}\}$ and $\tilde{r}= \sup_n \{ \tilde{t}- t_n\}$, then $\{x_n \}$ is contained in the subset
\[\bar B^+_{\bar{r}}(\bar{x})\cap
\bar B^-_{\tilde{r}}(\tilde{x}),\]
which is compact because it is the intersection of a compact subset and a closed subset.
Therefore we can extract a subsequence such that $(x_n,t_n)$ converges to $(x_0,t_0)$.
Now set $\bar{r}_n=\dist(\bar{x},x_n)$, $\tilde{r}_n=\dist(x_n,\tilde{x})$,  $\bar{r}_0=\dist(\bar{x},x_0)$  and $\tilde{r}_0=\dist(x_0,\tilde{x})$.
We know that  $t_n-\bar{t}\geq\bar{r}_n$ and $\tilde{t}- t_n\geq \tilde{r}_n$ and, as a consequence, we  have
$t_0-\bar{t}\geq\bar{r}_0$ and $\tilde{t}-t_0\geq \tilde{r}_0$. Hence it follows that

\[(x_0,t_0)\in \left(\bar B^+_{s_0(\tilde t-\bar t)} (\bar{x})\cap\bar B^-_{(1-s_0)(\tilde t-\bar t)} (\tilde{x})\right)\times\{\bar t+s_0(\tilde t-\bar t)\},\]
with $s_0=\frac{t_0-\bar t}{\tilde t-\bar t}$, and then it belongs to $J^+(p)\cap J^-(q)$.
Therefore $({\mathcal M},g)$ is globally hyperbolic.

 Now  we show (2).  We can assume without loss of generality that $\bar t=0$. We  will prove that $({\mathcal M}_0,F)$ is forward complete
 showing that every constant speed  geodesic $x\colon[0,b)\to {\mathcal M}_0$ can be extended to $b$.
 Assume that $x$ has  been parameterized with speed equal to $1$. Let $\{s_n\}\subset[0,b)$ be a sequence converging to $b$.
 We consider  the lightlike curve $\gamma:[0,b)\to {\mathcal M}$ such that $\gamma(s)=\big(x(s),-b+s\big)$.
 Then $ (x(\bar{s}),0)\in J^+(x(0),-b)$ for every $\bar{s}\in [0,b)$, because we can consider the lightlike curve $\gamma(s)$
 with $s\in [0,\bar s]$ and then the timelike curve $(x(\bar s), -b +s)$ with $s\in [\bar s, b]$.
Since in a globally hyperbolic manifold the intersection of the  future or the past of a point with a Cauchy surface is compact
(see for instance \cite[Proposition 6.6.6]{HawEll73}), the sequence  $x(s_n)$ is contained in a compact subset and converges
in contradiction with the fact that  $x$ is inextendible. Finally, arguing as above, we can show  that $({\mathcal M}_0,F)$ is also backward complete.
\qed\end{proof}

From Proposition~\ref{multiplicity}, we obtain the following result,
which gives a more geometrical interpretation  of previous results
\cite{FoGiMa95,Piccio97} because, apart  from the non-triviality of
the topology of the spacetime, it rests only on the completeness of
the Randers metric $F$. \bpr \label{lightrays} Let $(\mathcal M,g)$
be a conformally standard stationary Lorentzian manifold and
consider a point $(\bar x,\varrho_0)$ and the timelike curve
$\gamma(s)=(\tilde x,s)$. Assume that $(\mathcal M_0,F)$ is forward
or backward complete, then there exists a future-pointing light ray joining $(\bar
x,\varrho_0)$  and $\gamma(s)$. Moreover, assume that $\mathcal M_0$
is non-contractible, then there exist infinitely many lightlike
geodesics $z_n=(x_n, t_n)$ joining the point $(\bar x,\varrho_0)$
with the curve $\gamma(s)$ and having arrival time $T(x_n)\to
+\infty$, as $n\to\infty$. \epr \bere We observe that, since the
multiple geodesics found in the previous theorem have different
arrival time, they are also {\it geometrically distinct}. However we
cannot conclude, in general, that their spatial projections $x_n$
are geometrically distinct. \ere \bere Proposition~\ref{lightrays}
can be  generalized to lightlike geodesics joining two submanifolds
in $(\mathcal M,g)$ as in \cite{PerPic98}. Moreover,  as a closed geodesic exists on every compact Finsler manifold, we can also obtain
the existence of at least one non-trivial spatially periodic
lightlike geodesic, whenever $\mathcal M_0$ is compact. \ere 
\bere
A fully analogous result can be proved for past-pointing light rays by using the {\em reversed} Fermat metric 
\[F^*(x,y)=-\tig(x)[\delta(x),y]+
\sqrt{\tig(x)[\delta(x),y]^2+\tig(x)[y,y]}\]
and the arrival time functional $T^*(z)=t_0-\inte F^*(x,\dot x)\de s$. The reversed metric is related to the negative solution of Eq. \eqref{light}.  
We point out that multiplicity results about lightlike geodesics
connecting a point and a timelike curve in a spacetime are important
in the study of the gravitational lensing (see for instance
\cite{GiMaPi02,Perlic04}), that is, the deflection of light rays due
to the gravitational field of a galaxy. According to gravitational
lensing, the above result for past-pointing light rays can be interpreted as follows: $(\mathcal
M,g)$ is a conformally stationary spacetime having a  non trivial
topology, the point $(\bar x, \varrho_0)$ represents the position
and the time in which an observer receives the light signals, that
is, the lightlike geodesics emitted from a source whose trajectory in the spacetime is the curve $\gamma$. The fact that there exist infinitely many
lightlike geodesics connecting $(\bar x, \varrho_0)$ to $\gamma$
means that the observer sees, at the same instant of time,  many
images of the same source. \ere \bere\label{globalcondition} In view
of the importance of the Fermat metric completeness in the statement
of Proposition \ref{lightrays}, it is natural to ask under what
conditions on $g_0$, $\beta$ and $\delta$, the Fermat metric $F$ is
forward or backward complete. In the paper \cite{Sanche97}, it is
proved that a conformally standard stationary spacetime is  globally
hyperbolic, with Cauchy surface $\mathcal M_0\times\{0\}$, and then
by Theorem~\ref{gh} its Fermat metric is forward and backward 
complete, if $g_0$ is complete and $\beta$ and $|\delta|_0^2$ have
at most quadratic growth at infinity i.e. there exist constants
$c_1,\ c_2,\ c_3,\ c_4\geq 0$ such that \beq\label{optimal}
\begin{split}
&|\delta(x)|^2_0\leq c_1\disto^2(x,x_0)+c_2,\\
&\beta(x)\leq c_3\disto^2(x,x_0)+c_4,
\end{split}
\eeq
where $x_0$ is any fixed point in $\mathcal M_0$, $|\cdot|_0$ is the norm associated  to the metric $g_0$ and \disto is the distance on $\mathcal M_0$ induced by the metric $g_0$.
On the other hand, we can obtain a condition for the Fermat metric completeness directly from  Remark~\ref{completo}.
In fact,  it is enough to show that $g_0/\beta$ is complete and $\|\omega\|<1$.
Using the Cauchy-Schwarz inequality $g_0(y,y)\geq g_0(\delta,y)^2/ |\delta|^2_0$, we obtain a sufficient condition for $\|\omega\|<1$ as
\begin{equation}\label{betadelta}
\sup_{x\in\mathcal M_0}\frac{|\delta(x)|_0}{\sqrt{|\delta(x)|_0^2+\beta(x)}}<1.
\end{equation}
\end{remark}
\subsection{Timelike geodesics with fixed energy in stationary spacetimes}
In this subsection we reconsider   the Fermat  metric on a
one-dimensional  higher  manifold  in order to prove  multiplicity
of timelike geodesics with fixed energy  on  a standard
stationary spacetime $({\mathcal M},g)$, where $g$ is given by
\eqref{l}. Observe that, as timelike geodesics are not invariant under
conformal changes of the metric,  we are now obliged to consider only
standard stationary spacetimes. The idea is to use  
a Kaluza-Klein model  without the electromagnetic field (see  \cite{CapMas04} for an existence result of solutions for the relativistic Lorentz force equation based on Kaluza-Klein). More precisely, we seek for timelike geodesics $z\colon
[0,1] \rightarrow{\mathcal M}$ connecting a point $(x_0,t_0)\in
\mathcal M$ with a timelike curve $\gamma(s)=(x_1,s)\colon \R\to
\mathcal M$ and having a priori fixed energy $E_z = g(z)[\dot z,\dot
z]=-E<0$, for all $s\in[0,1]$.

We extend the Riemannian manifold ${\mathcal M}_0$ to the manifold ${\mathcal N}_0={\mathcal M}_0\times\R$ endowed with the metric
$n_0=g_0+\de u^2$, where  $u$ is the natural coordinate on $\R$,  and we associate
to the manifold ${\mathcal N}_0$ the one dimensional higher  Lorentzian manifold  $({\mathcal  N}, n)$, with the metric $n$ defined as
\begin{equation}\label{extend} n(x,u,t)[(y,v,\tau),(y,v,\tau)]=g_0(x)[y,y]+v^2 + 2g_0(x)[\delta(x),y]\tau -\beta(x)\tau^2.
\end{equation}
Since  $\partial_u$ is a Killing vector field for the metric  $n$, geodesics $\varsigma(s)=(x(s),u(s),t(s))$ in $({\mathcal N},n)$
have to satisfy also the conservation law 
\[n[\dot \varsigma,\partial_u]=\mathrm{const.},\]
which  implies that the $u$ component of  a  geodesic
is an affine function.
Moreover the projection $z(s) = (x(s),t(s))$ on ${\mathcal M}$ of  $\varsigma$ is a geodesic for $({\mathcal M},g)$.
In particular lightlike geodesics for the metric $n$ satisfy the following equation
\[g_0[\dot x,\dot x]+ 2g_0[\delta,\dot x]\dot t -\beta\dot t^2=-\dot u^2=\mathrm{const.}\]
Thus  in order to find timelike geodesics $z=(x,t)$ in $({\mathcal M},g)$ with fixed energy $-E<0$ it is enough to find
lightlike geodesics in $({\mathcal N},n)$ whose $u$ component has derivative equal to $\sqrt{E}$.
Fermat's principle in Subsection~\ref{lightlike} can be restated in
$({\mathcal N},n)$, reducing lightlike geodesics on $({\mathcal N},n)$ to geodesics for the Fermat  metric $\tilde F$ on the manifold ${\mathcal N}_0$, where $\tilde F$ is given by
\beq\label{RandersforE}
\tilde{F}((x,u),(y,v))=\sqrt{\tig[y,y]+v^2/\beta(x)+\tig[\delta(x),y]^2 }+\tig[\delta(x),y] ,
\eeq
for all $((x,u),(y,v))\in T{\mathcal N}_0$. We recall that $\tig=g_0/\beta$.
Therefore for any  value of energy $-E<0$ we obtain the following result,  which improves previous results about timelike geodesics with fixed energy on standard stationary Lorentzian manifolds as in  \cite{BaGeS02}, where $\delta=0$, and \cite{Germin06}, where only some ranges of values for $E$ are allowed.
\bpr\label{cauchypropertime}
Let $({\mathcal M},g)$ be a  standard stationary Lorentzian manifold. Assume that  $({\mathcal M}_0,F)$
is forward or backward complete and moreover assume that ${\mathcal M}_0$ is non-contractible, then
there exist infinitely many timelike geodesics $z_n=(x_n,t_n)$ connecting the point $(\bar x,t_0)\in \mathcal M$ with the timelike curve $\gamma(s)=(\tilde x,s)$,
parameterized on the interval $[a,b]$, having fixed energy $-E$ and  diverging arrival time.
\epr
\begin{proof}
Observe that if $\{(x_n,u_n)\}\subset {\mathcal N}_0$ is a forward Cauchy sequence for the Randers metric $\tilde F$ defined at \eqref{RandersforE},
then $\{x_n\}\subset M$  is a forward Cauchy sequence for the Fermat metric $F$ on $\mathcal M_0$ defined in \ref{Fermatmetric}.
Hence  $x_n$ converges and $\beta$ is bounded on $\{x_n\}$. Thus also $u_n$ is a Cauchy sequence in \R and therefore
$\{(x_n,u_n)\}$ converges, i.e. $({\mathcal N}_0,\tilde F)$ is forward complete. Then apply Proposition~\ref{multiplicity} to the Randers manifold
$({\mathcal N}_0,\tilde F)$
and to the functional $J((x,u))=\int_a^b \tilde F^2\big((x,u),(\dot x,\dot u)\big)\de s$ defined on the manifold \[\Lambda_{\{(\bar x,aE^{1/2})\}\times\{(\tilde x,bE^{1/2})\}}(\mathcal N_0)\]
(here the curves are parametrized on $[a,b]$) and use Fermat's principle on the manifold $({\mathcal N},n)$, between the point $(\bar x,aE^{1/2},t_0)$ and the curve
$s\mapsto(\tilde x,bE^{1/2},s)$.
\qed\end{proof}
The case  $E=1$  is  the most   interesting one, because timelike
geodesics with $E=1$ correspond to test particles, freely falling in
the gravitational field $g$,  parameterized with respect to the {\em
proper time} (see \cite{HawEll73}). In such a case, fixing the
interval of parameterization  is equivalent to fixing the arrival
proper time of the trajectory.

% BibTeX users please use
% \bibliographystyle{}
% \bibliography{}

\end{document}